# ANALYSIS OF VARIANCE, COEFFICIENT OF DETERMINATION AND $F$-TEST FOR LOCAL POLYNOMIAL REGRESSION

By Li-Shan Huang [1] and Jianwei Chen

*University of Rochester and San Diego State University*

This paper provides ANOVA inference for nonparametric local polynomial regression (LPR) in analogy with ANOVA tools for the classical linear regression model. A surprisingly simple and exact local ANOVA decomposition is established, and a local R-squared quantity is defined to measure the proportion of local variation explained by fitting LPR. A global ANOVA decomposition is obtained by integrating local counterparts, and a global R-squared and a symmetric projection matrix are defined. We show that the proposed projection matrix is asymptotically idempotent and asymptotically orthogonal to its complement, naturally leading to an $F$-test for testing for no effect. A by-product result is that the asymptotic bias of the "projected" response based on local linear regression is of quartic order of the bandwidth. Numerical results illustrate the behaviors of the proposed R-squared and $F$-test. The ANOVA methodology is also extended to varying coefficient models.

**1. Introduction.** Nonparametric regression methods such as local polynomial regression (LPR) (Fan and Gijbels [9], Wand and Jones [26]), smoothing splines (Eubank [8]) and penalized splines (Ruppert, Wand and Carroll [23]) are widely used to explore unknown trends in data analysis. Given the popularity of these methods, a set of analysis of variance (ANOVA) inference tools, analogous to those of linear models, will be very useful in providing interpretability for nonparametric curves. In this paper, we aim to develop ANOVA inference for LPR. Some of the work in this paper was motivated by the authors' consulting project experiences, where clients presented with a nonparametric smooth curve would frequently ask if there would be an ANOVA table explicitly summarizing the fitted curve by sums of squares,

Received March 2005; revised July 2007.
[1]Supported in part by DOD MURI Grant FA9550-04-1-0430.

degrees of freedom, and an $F$-test for no effect. In addition, we are interested in exploring a geometric representation and establishing a projection view for LPR.

Consider a simple bivariate case: data $(X_i, Y_i)$, $i = 1, \ldots, n$, are drawn from the model

$$Y = m(X) + \sigma(X)\varepsilon, \tag{1.1}$$

where $X$ and $\varepsilon$ are independent, and $\varepsilon$ has a mean 0 and unit variance. Typically one is interested in estimating the conditional mean, $m(x) = E(Y|X = x)$, while the conditional variance is $\sigma^2(x) = \text{Var}(Y|X = x)$. The theoretical ANOVA decomposition for (1.1) is

$$\text{Var}(Y) = \int (m(x) - \mu_y)^2 f(x)\,dx + \int \sigma^2(x) f(x)\,dx, \tag{1.2}$$

where $f(x)$ is the underlying density function for $X_1, \ldots, X_n$, and $\mu_y$ denotes the unconditional expected value of $Y$. Below we review briefly some related work on ANOVA inference for nonparametric regression.

In LPR literature, we are not aware of a sample ANOVA decomposition for (1.2). A commonly used residual sum of squares (RSS) is $\sum_{i=1}^{n}(Y_i - \hat{m}(X_i))^2$, where $\hat{m}(X_i)$ denotes a nonparametric estimate for $m(X_i)$, $i = 1, \ldots, n$, but RSS is not associated with a valid ANOVA decomposition, in the sense that generally $\sum_{i=1}^{n}(Y_i - \bar{Y})^2 \neq \sum_{i=1}^{n}(\hat{m}(X_i) - \bar{Y})^2 + \sum_{i=1}^{n}(Y_i - \hat{m}(X_i))^2$, where $\bar{Y}$ is the sample mean of $Y_i$'s. Ramil-Novo and González-Manteiga [22] established an ANOVA decomposition for smoothing splines with a bias term. An ANOVA-related quantity is the R-squared, or the coefficient of determination. Theoretically, it measures $\eta^2 = 1 - E(\text{Var}(Y|X))/\text{Var}(Y) = \text{Var}(E(Y|X))/\text{Var}(Y)$. Doksum and Samarov [5] suggested an estimate

$$R_\rho^2 = \frac{[n^{-1}\sum_i(\hat{m}(X_i) - \bar{m})(Y_i - \bar{Y})]^2}{[n^{-1}\sum_i(\hat{m}(X_i) - \bar{m})^2][n^{-1}\sum_i(Y_i - \bar{Y})^2]}, \tag{1.3}$$

where $\bar{m} = n^{-1}\sum_i \hat{m}(X_i)$. However, the correlation-based $R_\rho^2$ does not possess an ANOVA structure. For a local version of the R-squared measure, see Bjerve and Doksum [3], Doksum et al. [4] and Doksum and Froda [6]. An attempt to provide an analogous projection matrix is the so-called "smoother matrix" $\mathbf{S}$, $n \times n$, so that $\mathbf{Sy} = \hat{\mathbf{m}}$ with $\mathbf{y} = (Y_1, \ldots, Y_n)^T$ and $\hat{\mathbf{m}} = (\hat{m}(X_1), \ldots, \hat{m}(X_n))^T$. See, for example, Hastie and Tibshirani [13]. However, $\mathbf{S}$ lacks for properties of a projection matrix; it is non-idempotent and nonsymmetric in the case of local linear regression. Another essential ANOVA element is the degree of freedom (DF). Hastie and Tibshirani [13] discussed three versions: $\text{tr}(\mathbf{S})$, $\text{tr}(\mathbf{S}^T\mathbf{S})$ and $\text{tr}(2\mathbf{S} - \mathbf{S}^T\mathbf{S})$, where "tr" denotes the trace of a matrix. Zhang [27] gave asymptotic expressions on DF for LPR. On testing for no effect, Azzalini, Bowman and Hardle [1], Hastie and Tibshirani [13] and Azzalini and Bowman [2] introduced tests with the $F$-type form of test statistics based on RSS. Fan, Zhang and Zhang [10] established the generalized likelihood ratio test with an $F$-type test statistic and an asymptotic chi-square distribution. Other $F$-flavor tests include Gijbels and Rousson [12].

From the discussion above, we believe that there is a need to further investigate an ANOVA framework for LPR. Our focus on LPR arises naturally since it is a "local" least squares technique. A surprisingly simple local ANOVA decomposition is established in Section 2, leading naturally to defining a local R-squared. Then by integrating local counterparts, a global ANOVA decomposition is established, from which a global R-squared and a symmetric matrix $H^*$, like a projection matrix, are defined. We note that the proposed global SSE (sum of squares due to error) is the same as the "smooth backfitting" error given in Mammen, Linton and Nielsen [19] and Nielsen and Sperlich [20] for estimation under generalized additive models (Hastie and Tibshirani [13]). We show that when conditioned on $\{X_1, \ldots, X_n\}$, $H^*$ is asymptotically idempotent and $H^*$ and its complement $(I - H^*)$ are asymptotically orthogonal, leading naturally to an $F$-test for testing no effect. A by-product is that the conditional bias of the "projected" response $H^*\mathbf{y}$ based on local linear regression is of order $h^4$, with $h$ the bandwidth. To show that the ANOVA framework can be extended to the multivariate case, expressions of local and global ANOVA decomposition are derived for varying coefficient models (VCM) (Hastie and Tibshirani [14]) in Section 3. Section 4 contains numerical results on the performance of the proposed global R-squared and the $F$-test for no effect. In summary, our results are under one framework containing all essential ANOVA elements: (i) a local exact ANOVA decomposition, (ii) a local R-squared, (iii) a global ANOVA decomposition, (iv) a global R-squared, (v) an asymptotic projection matrix $H^*$, (vi) nonparametric degree of freedom defined by $\text{tr}(H^*)$ and (vii) an $F$-test for testing no effect. The results also give new insights of LPR being a "calculus" extension of classical polynomial models and provide a new geometric view on LPR highlighted by $H^*$. Extension of the ANOVA inference to partially linear models, generalized additive models and semiparametric models is in progress.

**2. ANOVA for local polynomial regression.** We begin by introducing LPR (Fan and Gijbels [9], Wand and Jones [26]) under (1.1). Assume that locally for data $X_i$'s in a neighborhood of $x$, $m(X_i)$ can be approximated by $m(x) + m'(x)(X_i - x) + \cdots + m^{(p)}(x)(X_i - x)^p/p!$, based on a Taylor expansion. Then this local trend is fitted by weighted least squares as the following:

$$(2.1) \qquad \min_{\beta} n^{-1} \sum_{i=1}^{n} \left(Y_i - \sum_{j=0}^{p} \beta_j(X_i - x)^j\right)^2 K_h(X_i - x),$$

where $\beta = (\beta_0, \ldots, \beta_p)^T$, $K_h(\cdot) = K(\cdot/h)/h$, and the dependence of $\beta$ on $x$ and $h$ is suppressed. The function $K(\cdot)$ is a nonnegative weight function, typically a symmetric probability density function, and $h$ is the smoothing parameter, determining the neighborhood size for local fitting. Let $\hat{\beta} = (\hat{\beta}_0, \ldots, \hat{\beta}_p)^T$ denote the solution to (2.1). It is clear that $\hat{\beta}_0$ estimates $m(x)$ of interest and $j!\hat{\beta}_j$ estimates the $j$th derivative $m^{(j)}(x)$, $j = 1, \ldots, p$. For convenience of developing ANOVA inference in this paper, we define a local SSE as the resulting (2.1) divided by the sum of local weights:

$$(2.2) \quad SSE_p(x;h) = \frac{n^{-1}\sum_{i=1}^n (Y_i - \sum_{j=0}^p \hat{\beta}_j(X_i - x)^j)^2 K_h(X_i - x)}{n^{-1}\sum_{i=1}^n K_h(X_i - x)}.$$

The denominator of (2.2) is the kernel density estimator $\hat{f}(x;h)$ (Silverman [24]) for $f(x)$. Similar treatment can be found in Qiu [21], so that $SSE_p(x;h)$ estimates $\sigma^2(x)$. We note that (2.2) is equivalent to the SSE for weighted least squares regression given in Draper and Smith [7].

Recall that in the linear regression setting, the sample ANOVA decomposition is given as $SST \equiv n^{-1}\sum_i (Y_i - \bar{Y})^2 = n^{-1}\sum_i (\hat{Y}_i - \bar{Y})^2 + n^{-1}\sum_i (Y_i - \hat{Y}_i)^2 \equiv SSR + SSE$, where $\hat{Y}_i$'s denote fitted values for $Y_i$'s from a linear model, $SST$ the corrected sum of squares for $Y_i$'s, and $SSR$ the sum of squares due to regression. In the literature of weighted least squares regression (e.g., Draper and Smith [7]) with weight $w_i$ assigned to $(X_i, Y_i)$, the sample ANOVA decomposition is

$$(2.3) \quad \sum_i (Y_i - \bar{Y}_w)^2 w_i = \sum_i (\hat{Y}_{i,w} - \bar{Y}_w)^2 w_i + \sum_i (Y_i - \hat{Y}_{i,w})^2 w_i,$$

where $\bar{Y}_w = \sum_i Y_i w_i / \sum_i w_i$ and $\hat{Y}_{i,w}$ is the resulting fitted value for $Y_i$.

2.1. *Local ANOVA decomposition and a pointwise R-squared.* The local least squares feature of LPR leads us to consider whether an analogous local (pointwise) ANOVA decomposition exists. We note that it is not suitable to adopt (2.3) directly. By forcing a local fit of $\bar{Y}$, we obtain a finite-sample and exact local ANOVA decomposition in Theorem 1 for LPR. In addition to $SSE_p(x;h)$ in (2.2), local SST and local SSR are defined as follows:

$$(2.4) \quad \begin{aligned} SST(x;h) &= \frac{n^{-1}\sum_{i=1}^n (Y_i - \bar{Y})^2 K_h(X_i - x)}{\hat{f}(x;h)}, \\ SSR_p(x;h) &= \frac{n^{-1}\sum_{i=1}^n (\sum_{j=0}^p \hat{\beta}_j(X_i - x)^j - \bar{Y})^2 K_h(X_i - x)}{\hat{f}(x;h)}. \end{aligned}$$

Note that both $SSE_p(x;h)$ and $SSR_p(x;h)$ use all the fitted parameters $\hat{\beta}_j$'s, in contrast to RSS using only $\hat{\beta}_0$.

THEOREM 1. *An exact and finite-sample ANOVA decomposition is obtained for local polynomial fitting at a grid point $x$ in the range of $X_i$'s:*

$$(2.5) \qquad SST(x;h) = SSE_p(x;h) + SSR_p(x;h).$$

*In addition, $SSE_1(x;h)$ for local linear regression ($p=1$) is related to the weighted least squared error of the Nadaraya–Watson estimator ($p=0$), as given below:*

$$(2.6) \qquad \begin{aligned} SSE_1(x;h) &= \frac{n^{-1}\sum(Y_i - \hat{m}_{NW}(x))^2 K_h(X_i - x)}{\hat{f}(x;h)} \\ &\quad - \hat{\beta}_1^2 \frac{n^{-1}\sum(X_i - \bar{X}_k)^2 K_h(X_i - x)}{\hat{f}(x;h)}, \end{aligned}$$

*where $\hat{m}_{NW}(x) = n^{-1}\sum_i K_h(X_i - x)Y_i/\hat{f}(x;h)$ and $\bar{X}_k = n^{-1}\sum_i X_i K_h(X_i - x)/\hat{f}(x;h)$.*

The proof of Theorem 1 is mainly algebraic and hence is omitted; (2.6) is simply (2.3). The "exact" expression (2.5) is very attractive and has an appealing interpretation of comparing the local fit with the simple no-effect $\bar{Y}$ in the same local scale. It is easy to see that $SSR_p(x;h)$ estimates $(m(x) - \mu_y)^2$ and $SSE_p(x;h)$ estimates $\sigma^2(x)$.

Based on (2.5), we define a local (pointwise) R-squared at $x$ as follows:

$$(2.7) \qquad R_p^2(x;h) = 1 - \frac{SSE_p(x;h)}{SST(x;h)} = \frac{SSR_p(x;h)}{SST(x;h)}.$$

From Theorem 1, $R_p^2(x;h)$ is always between 0 and 1, and $R_1^2(x;h)$ for local linear regression is always greater than $R_0^2(x;h)$ for the Nadaraya–Watson estimator with the same bandwidth and kernel function. A plot of $R_p^2(x;h)$ versus $x$ will give an idea of the quality of estimation at different regions of data. $R_p^2(x;h)$ is a measure for the proportion of local variation explained by the local polynomial fit. We note that $R_p^2(x;h)$ is invariant with respect to linear transformations of $Y_i$'s, and will be invariant for linear transformations of $X_i$'s ($aX_i + b$) if the bandwidth is taken proportional to the transformation, $ah$, accordingly. The classical R-squared for polynomial models can be viewed as a special case of (2.7), when using the uniform kernel at only one grid point $\bar{X}$. Thus LPR, fitting local polynomials across data, is like a calculus extension of classical polynomial models.

2.2. *Global ANOVA decomposition and coefficient of determination.* We now turn to developing a global ANOVA decomposition. It is convenient to introduce some conditions here.

CONDITIONS (A).

(A1) The design density $f(x)$ is bounded away from 0 and $\infty$, and $f(x)$ has a continuous second derivative on a compact support.
(A2) The kernel $K(\cdot)$ is a Lipschitz continuous, bounded and symmetric probability density function, having a support on a compact interval, say $[-1, 1]$.
(A3) The error $\varepsilon$ is from a symmetric distribution with mean 0, variance 1, and a finite fourth moment.
(A4) The $(p+1)$st derivative of $m(\cdot)$ exists.
(A5) The conditional variance $\sigma^2(\cdot)$ is bounded and continuous.

Based on (2.5), a global ANOVA decomposition can be established by integrating local counterparts in (2.5):

$$
\begin{aligned}
SST(h) &= \int SST(x;h)\hat{f}(x;h)\,dx, \\
SSE_p(h) &= \int SSE_p(x;h)\hat{f}(x;h)\,dx, \\
SSR_p(h) &= \int SSR_p(x;h)\hat{f}(x;h)\,dx.
\end{aligned}
\tag{2.8}
$$

Then a global ANOVA decomposition is

$$SST = SSE_p(h) + SSR_p(h), \tag{2.9}$$

which corresponds to the theoretical version (1.2). Since $\int K_h(X_i - x)\,dx = 1$ under Conditions (A1) and (A2), $\int SST(x;h)\hat{f}(x;h)\,dx = n^{-1}\sum_{i=1}^{n}(Y_i - \bar{Y})^2 = SST$ in (2.9). We then define a global R-squared as

$$R_p^2(h) = 1 - \frac{SSE_p(h)}{SST} = \frac{SSR_p(h)}{SST}, \tag{2.10}$$

and we name it the "ANOVA" R-squared. We further investigate some asymptotic properties of $R_p^2(h)$. For simplicity, we focus on the case of an odd degree, for example, $p = 1$, in Theorem 2. A by-product is that $SSE(h)$ is a $\sqrt{n}$-consistent estimate for $\sigma^2$ when assuming homoscedasticity.

THEOREM 2. *Assume that as $n \to \infty$, $h = h(n) \to 0$. When fitting LPR with an odd $p$, under Conditions (A) with $nh^{2p+2} \to 0$ and $nh^2 \to \infty$:*

(a) *The asymptotic conditional bias of $R_p^2(h)$ is*

$$-h^2 \frac{\mu_2}{2\sigma_y^2} \int \sigma^2(x) f''(x)\,dx (1 + o_P(1)).$$

(b) *The asymptotic conditional variance of $R_p^2(h)$ is*

$$n^{-1}\left(\frac{\mathrm{Var}(\varepsilon^2)}{\sigma_y^4}E(\sigma^4(X))\left(\int K_0^*(v)\,dv\right) + \frac{(m_4 - \sigma_y^4)(E(\sigma^2(X)))^2}{\sigma_y^8}\right)$$
$$\times (1 + o_P(1)),$$

where $\sigma_y^2$ is the variance of $Y$, $m_4 = E\{(Y - \mu_y)^4\}$ is the fourth central moment of $Y$, and $K_0^*(v) = \int K(u)K(v-u)\,du$ denotes the convolution of $K$ and itself.

(c) *Under the homoscedastic assumption and conditioned on* $\{X_1, \ldots, X_n\}$, $R_p^2(h)$ *converges in distribution to a normal distribution with the above asymptotic conditional bias and variance.*

(d) *Under the assumptions in* (c), $SSE_p(h)$ *is a $\sqrt{n}$-consistent estimate for $\sigma^2$. Its asymptotic conditional bias is $o_P(n^{-1/2})$ if $\int f''(x)\,dx = 0$ and its asymptotic conditional variance $n^{-1}\sigma^4(\int K_0^*(v)\,dv)(1 + o_P(1))$.*

Theorem 2 is a special case of Theorem 6 in Section 3, and hence the proof of Theorem 6 (in the Appendix) is applicable to Theorem 2. The condition on the bandwidth in Theorem 2 becomes $h = o(n^{-1/4})$ and $n^{-1/2} = o(h)$ for the case of $p = 1$. It is known that the optimal bandwidth for estimating $m(\cdot)$ with $p = 1$ is of order $n^{-1/5}$ (e.g., Fan and Gijbels [9]). It is not surprising that we need a smaller bandwidth than the rate of $n^{-1/5}$ to obtain a $\sqrt{n}$-consistent estimate for $\sigma^2$.

2.3. *Asymptotic projection matrix.* Under Conditions (A1) and (A2), $SSE_p(h)$ and $SSR_p(h)$ can be rewritten as

$$SSE_p(h) = n^{-1}\left\{\sum_i Y_i^2 - \int \sum_i \left(\sum_j \hat{\beta}_j(x)(X_i - x)^j\right)^2 K_h(X_i - x)\,dx\right\},$$

$$SSR_p(h) = n^{-1}\left\{\int \sum_i \left(\sum_j \hat{\beta}_j(x)(X_i - x)^j\right)^2 K_h(X_i - x)\,dx - n\bar{Y}^2\right\},$$

and in a matrix expression,

(2.11) $\quad SSE_p(h) = n^{-1}\mathbf{y}^T(I - H^*)\mathbf{y}, \qquad SSR_p(h) = n^{-1}\mathbf{y}^T(H^* - L)\mathbf{y},$

where $L$ is an $n \times n$ matrix with entries $1/n$. In this subsection, we further explore if $H^*$ behaves like a projection matrix. The $H^*$ matrix can be written as $H^* = \int WH\hat{f}(x;h)\,dx$, where $W$ is a diagonal matrix with entries $K_h(X_i - x)/\hat{f}(x;h)$, $H = \mathbf{X}(\mathbf{X}^T W\mathbf{X})^{-1}\mathbf{X}^T W$ is the local projection matrix for (2.1) with $\mathbf{X}$ the design matrix for (2.1), and the integration is performed element by element in the resulting matrix product. $H^*$ depends on the data points $X_i$'s, kernel function $K$ and the bandwidth $h$. Under Conditions (A1) and (A2), $H^*\mathbf{1} = \mathbf{1}$, where $\mathbf{1}$ denotes an $n$-vector of 1's. Therefore the projected response $H^*\mathbf{y} = \mathbf{y}^*$ is a vector with each element $Y_i^*$ being a weighted average of $Y_i$'s. The matrix $H^*$ is clearly a symmetric matrix, but it is not idempotent. Given this fact, we take a step back to explore if $H^*$ is asymptotically idempotent when conditioned on $\{X_1, \ldots, X_n\}$.

The authors are not aware of standard criteria of asymptotic idempotency. Below we define a criterion for asymptotic idempotency and asymptotic orthogonality in a nonparametric regression setting:

DEFINITION. 1. Conditioned on $\{X_1, \ldots, X_n\}$, an $n \times n$ matrix $A_n$ is asymptotically idempotent, if for any random $n$-vector response $\mathbf{y}$ with finite expected value, $E\{(A_n - A_n^2)\mathbf{y}|X_1, \ldots, X_n\}$ tends to a zero vector in probability as $n \to \infty$, that is, each element of $E\{(A_n - A_n^2)\mathbf{y}|X_1, \ldots, X_n\}$ is asymptotically zero in probability as $n \to \infty$.

2. Conditioned on $\{X_1, \ldots, X_n\}$, for two $n \times n$ matrices $A_n$ and $B_n$, they are asymptotically orthogonal, if for any random $n$-vector response $\mathbf{y}$ with finite expected value, $E\{A_n B_n \mathbf{y}|X_1, \ldots, X_n\}$ tends to a zero vector in probability as $n \to \infty$, that is, each element of $E\{A_n B_n \mathbf{y}|X_1, \ldots, X_n\}$ is asymptotically zero in probability as $n \to \infty$.

Denote the multiplier for $h^{p+1}$ ($p$ odd) or $h^{p+2}$ ($p$ even) in the first-order term for the conditional bias of $\hat{\beta}_0(x; h, p)$ as $b_{0,p}(x)$ (see Wand and Jones [26], page 125). The following theorem gives the rate of each element in $(H^* - H^{*2})\mathbf{y}$.

THEOREM 3. *Under Conditions (A), suppose local polynomial regression of order $p$ is fitted to data. The bandwidth $h \to 0$ and $nh \to \infty$, as $n \to \infty$.*

(a) *For $p \neq 1$, the asymptotic conditional bias of $Y_i^*$, $E\{Y_i^* - m(X_i)|X_1, \ldots, X_n\}$, for $i = 1, \ldots, n$, is at most*

$$(2.12) \quad \begin{cases} O(h^{p+1})(1 + o_P(1)), & p \text{ is odd}; \\ O(h^{p+2})(1 + o_P(1)), & p \text{ is even}. \end{cases}$$

(b) *For $p = 1$, the asymptotic conditional bias of $Y_i^*$, $i = 1, \ldots, n$, is of order $h^4$; more explicitly*

$$(2.13) \quad h^4\left(\frac{\mu_2^2 - \mu_4}{4}\right)\left\{m^{(4)}(X_i) + 2m^{(3)}(X_i)\frac{f'(X_i)}{f(X_i)} + m^{(2)}(X_i)\frac{f''(X_i)}{f(X_i)}\right\} + o_P(h^4).$$

(c) *Each element in $E\{(H^* - H^{*2})\mathbf{y}|X_1, \ldots, X_n\}$ is at most of order $O(h^{p+1}) \times (1 + o_P(1))$ for an odd $p$ with $p \geq 3$, at most $O(h^{p+2})(1 + o_P(1))$ if $p$ is even, and $O(h^4)(1 + o_P(1))$ when $p = 1$. Thus, conditioned on $\{X_1, \ldots, X_n\}$, $H^*$ is asymptotically idempotent and asymptotically a projection matrix.*

(d) *For local linear regression, the asymptotic conditional variance of $Y_i^*$ retains the order of $n^{-1}h^{-1}$:*

(2.14) $\quad \operatorname{Var}\{Y_i^* | X_1, \ldots, X_n\} = n^{-1}h^{-1}(1 + o_P(1))\kappa_0 \sigma^2 / f(X_i),$

*where $\kappa_0 = \int K_0^{*2}(v)\,dv - \frac{2}{\mu_2} \int K_0^*(v) K_1^*(v)\,dv + \frac{1}{\mu_2^2} \int K_1^{*2}(v)\,dv$ with $K_1^*(\cdot)$ the convolution of $uK(u)$ and itself.*

Theorem 3(b) implies that using the matrix $H^*$, one can achieve a surprising bias reduction effect for local linear regression, from the order of $h^2$ to $h^4$. While achieving bias reduction, the asymptotic conditional variance of $Y_i^*$ increases in the case of local linear regression. We calculate the constant term in (2.14) for the Epanechnikov and Gaussian kernel functions, and the ratios of the constant factors of $Y_i^*$ and local linear estimator $\hat{\beta}_0(X_i)$ are 1.38 and 1.10, respectively. It is of interest to know the form of $\mathbf{y}^*$,

$$(2.15) \quad H^*\mathbf{y} = \begin{pmatrix} \int (\hat{\beta}_0(x) + \cdots + \hat{\beta}_p(x)(X_1 - x)^p) K_h(X_1 - x)\,dx \\ \vdots \\ \int (\hat{\beta}_0(x) + \cdots + \hat{\beta}_p(x)(X_n - x)^p) K_h(X_n - x)\,dx \end{pmatrix}.$$

The projection $H^*\mathbf{y}$ uses all the fitted $\hat{\beta}_j(x)$'s through integration and the gain is reduction in the asymptotic bias. It is in contrast with $\hat{\beta}_0(X_i)$, which fits local polynomial at $X_i$ and throws away other fitted parameters when $p \geq 1$.

2.4. *An F-test for testing no effect.* Results in Theorem 3 naturally lead us to consider an $F$-test for testing no effect. The next theorem proposes an $F$-test that inherits properties of the classical $F$-tests.

THEOREM 4. *Under the conditions in Theorem 3 and conditioned on $\{X_1, \ldots, X_n\}$:*

(a) $(I - H^*)$ *and* $(H^* - L)$ *are asymptotically orthogonal, in the sense that*

$$E\{(I - H^*)(H^* - L)\mathbf{y} | X_1, \ldots, X_n\} = E\{(H^* - H^{*2})\mathbf{y} | X_1, \ldots, X_n\},$$

*which tends to a zero vector in probability.*

(b) *Under the simple homoscedastic assumption, an F-statistic is formed as*

$$(2.16) \quad F = \frac{SSR_p / (\operatorname{tr}(H^*) - 1)}{SSE_p / (n - \operatorname{tr}(H^*))},$$

*where $\operatorname{tr}(H^*)$ is the trace of $H^*$. Conditioned on $\{X_1, \ldots, X_n\}$, with the normal error assumption, the F-statistic (2.16) is asymptotically F-distributed with degrees of freedom $(\operatorname{tr}(H^*) - 1, n - \operatorname{tr}(H^*))$.*

TABLE 1
*ANOVA table for local polynomial regression*

| Source | Degree of freedom | Sum of squares | Mean squares | $F$ |
|---|---|---|---|---|
| Regression | $\operatorname{tr}(H^*) - 1$ | $SSR_p = n^{-1}\mathbf{y}^T(H^* - L)\mathbf{y}$ | $MSR_p = \frac{nSSR_p}{(\operatorname{tr}(H^*)-1)}$ | $\frac{MSR_p}{MSE_p}$ |
| Residual | $n - \operatorname{tr}(H^*)$ | $SSE_p = n^{-1}\mathbf{y}^T(I - H^*)\mathbf{y}$ | $MSE_p = \frac{n \times SSE_p}{(n-\operatorname{tr}(H^*))}$ | |
| Total | $(n-1)$ | $SST = n^{-1}\mathbf{y}^T(I - L)\mathbf{y}$ | | |

(c) *The conditional trace of $H^*$ for local linear regression is asymptotically*

$$(2.17) \quad \operatorname{tr}(H^*) = h^{-1}|\Omega|(\nu_0 + \nu_2/\mu_2)(1 + o_P(1)),$$

*where $|\Omega|$ denotes the range of $X_i$'s and $\nu_j = \int u^j K^2(u)\,du$.*

We remark that when a local $p$th polynomial approximation is exact, $E\{(H^* - H^{*^2})\mathbf{y}|X_1,\ldots,X_n\} = 0$, that is, $H^*$ is idempotent and the resulting $F$-statistic (2.16) has an exact $F$-distribution as in the classical settings. Based on (2.8) and Theorems 3 and 4, an ANOVA table for LPR is given in Table 1. It has been shown in Theorem 2(d) that $SSE_p(h)$ is a $\sqrt{n}$-consistent estimate for $\sigma^2$ when the error variance is homoscedastic. Table 1 shows that $MSE_p(h) = SSE_p(h)\frac{n}{n-\operatorname{tr}(H^*)}$ is an unbiased estimate for $\sigma^2$ in finite-sample settings, which is similar to the classical *MSE* in linear models. With the ANOVA table, an analogous adjusted R-squared may be defined as

$$(2.18) \quad R^2_{p,adj}(h) = 1 - \frac{SSE_p(h)/(n - \operatorname{tr}(H^*))}{SST/(n-1)}.$$

**3. Extension to varying coefficient models.** In this section, we extend the ANOVA decomposition to VCM, illustrating that the ANOVA framework can be extended to the multivariate case. Though there is no room in this paper for a full discussion of VCM, we develop expressions for local and global ANOVA decomposition and the ANOVA R-squared in this section.

The VCM assumes the following conditional linear structure:

$$(3.1) \quad Y = \sum_{k=1}^{d} a_k(U)X_k + \sigma(U)\varepsilon,$$

where $X_1,\ldots,X_d$, $d \geq 1$, are the covariates with $X_1 = 1$, $\mathbf{a}(U) = (a_1(U),\ldots,a_d(U))^T$ is the functional coefficient vector, $U$ and $\varepsilon$ are independent, and $\varepsilon$ has a mean 0 and unit variance. Specifically, when $d = 1$, model (3.1) is reduced to the bivariate nonparametric model (1.1). On the other hand, if the varying coefficients are constants, that is, $a_k(U) = a_k$, $k = 1,\ldots,d$, the model is the multivariate linear model. Based on (3.1), the theoretical ANOVA decomposition is

$$\text{Var}(Y) = \text{Var}(E(Y|U, X_1, \ldots, X_d)) + E(\text{Var}(Y|U, X_1, \ldots, X_d))$$
$$(3.2) \qquad = \int (\mathbf{a}(u)^T \mathbf{x} - \mu_y)^2 f(\mathbf{x}|u) g(u) \, d\mathbf{x} \, du + \int \sigma^2(u) g(u) \, du,$$

where $g(u)$ denotes the underlying density function for $U$, and $f(\mathbf{x}|u)$ the underlying conditional density function of $\mathbf{x} = (X_1, \ldots, X_d)^T$ given $u$.

Hoover et al. [15] and Fan and Zhang [11] applied LPR to estimate the varying-coefficient function vector $\mathbf{a}(U)$. Assume that the $(p+1)$st-order derivative of $\mathbf{a}(U)$ exists, and data $(U_i, X_{i1}, \ldots, X_{id}, Y_i)$, $i = 1, \ldots, n$, are drawn from model (3.1). Based on a Taylor expansion, $a_k(U_i)$, $i = 1, \ldots, n, k = 1, \ldots, d$, is approximated by $\beta_{k,0}(u) + \beta_{k,1}(u)(U_i - u) + \cdots + \beta_{k,p}(u)(U_i - u)^p$, for $U_i$ in a neighborhood of a grid point $u$. Then local polynomial estimator $\hat{\beta}_k = (\hat{\beta}_{k,0}, \ldots, \hat{\beta}_{k,p})^T$, $k = 1, \ldots, d$, for VCM can be obtained by the following locally weighted least squares equation:

$$(3.3) \quad \min_{\beta} n^{-1} \sum_{i=1}^{n} \left( Y_i - \sum_{k=1}^{d} \sum_{j=0}^{p} \beta_{k,j} (U_i - u)^j X_{ik} \right)^2 K_h(U_i - u) / \hat{g}(u; h),$$

where $\beta = (\beta_{1,0}, \ldots, \beta_{1,p}, \ldots, \beta_{d,0}, \ldots, \beta_{d,p})^T$, and $\hat{g}(u; h) = n^{-1} \times \sum_{i=1}^{n} K_h(U_i - u)$ denotes the kernel density estimate for $g(u)$. For convenience, (3.3) and its solution are expressed in a matrix form. Let

$$\mathbf{X}_u = \begin{pmatrix} X_{11} & \cdots & X_{11}(U_1 - u)^p & \cdots & X_{1d} & \cdots & X_{1d}(U_1 - u)^p \\ \vdots & \ddots & \vdots & \ddots & \vdots & \ddots & \vdots \\ X_{n1} & \cdots & X_{n1}(U_n - u)^p & \cdots & X_{nd} & \cdots & X_{nd}(U_n - u)^p \end{pmatrix}_{n \times (p+1)d},$$

and $W_u$ be an $n \times n$ diagonal matrix of weights with $i$th element $K_h(U_i - u)/\hat{g}(u; h)$. Then the solution to (3.3) can be expressed as $\hat{\boldsymbol{\beta}}(u) = (\mathbf{X}_u^T W_u \mathbf{X}_u)^{-1} \times \mathbf{X}_u^T W_u \mathbf{y}$, and the local polynomial estimator for $\mathbf{a}(u)$ is

$$\hat{\mathbf{a}}(u) = (I_d \otimes e_{(p+1),1})(\mathbf{X}_u^T W_u \mathbf{X}_u)^{-1} \mathbf{X}_u^T W_u \mathbf{y},$$

where $\otimes$ denotes the Kronecker product and $e_{(p+1),k}$ is a $(p+1)$-dimension vector with 1 on the $k$th position and 0 elsewhere, and $\hat{\mathbf{a}}(u) = (\hat{\beta}_{1,0}(u), \ldots, \hat{\beta}_{d,0}(u))^T$.

Similarly to the bivariate case, Theorem 5 gives the local finite-sample ANOVA decomposition for VCM.

THEOREM 5. *Under model (3.1), an exact and finite-sample ANOVA decomposition is obtained for local polynomial fitting at a grid point $u$:*

$$SST(u; h) \equiv \frac{n^{-1} \sum_{i=1}^{n} (Y_i - \bar{Y})^2 K_h(U_i - u)}{\hat{g}(u; h)}$$

$$= \frac{n^{-1}\sum_{i=1}^{n}(Y_i - \hat{Y}_i(u))^2 K_h(U_i - u)}{\hat{g}(u;h)}$$

$$+ \frac{n^{-1}\sum_{i=1}^{n}(\hat{Y}_i(u) - \bar{Y})^2 K_h(U_i - u)}{\hat{g}(u;h)}$$

$$\equiv SSE_p(u;h) + SSR_p(u;h),$$

where $\hat{Y}_i(u) = e_{ni}\mathbf{X}_u(\mathbf{X}_u^T W_u \mathbf{X}_u)^{-1}\mathbf{X}_u^T W_u \mathbf{y} = \sum_{k=1}^{d}\sum_{j=0}^{p}\hat{\beta}_{k,j}(U_i-u)^j X_{ik}$ with $e_{ni}$ an n-dimension vector with 1 at the ith position and 0 elsewhere.

The ANOVA decomposition in Theorem 5 extends the bivariate ANOVA decomposition (2.5) to VCM in a straightforward way. A global ANOVA decomposition can be constructed by integrating the local counterparts in Theorem 5:

(3.4) $$SST(h) = SSE_p(h) + SSR_p(h),$$

where

$$SST = \int SST(u;h)\hat{g}(u;h)\, du = n^{-1}\mathbf{y}^T(I - L)\mathbf{y},$$

(3.5) $$SSE_p(h) = \int SSE_p(u;h)\hat{g}(u;h)\, du = n^{-1}\mathbf{y}^T(I - H_u^*)\mathbf{y},$$

$$SSR_p(h) = \int SSE_p(u;h)\hat{g}(u;h)\, du = n^{-1}\mathbf{y}^T(H_u^* - L)\mathbf{y},$$

where $H_u^* = \int W_u H_u \hat{g}(u;h)\, du$ is a symmetric $n \times n$ matrix with $H_u = \mathbf{X}_u(\mathbf{X}_u^T \times W_u \mathbf{X}_u)^{-1}\mathbf{X}_u^T W_u$. The matrix expression in the right-hand side of (3.5) is derived under Conditions (B1) below and (A2), and similarly to Section 2, $SST$ is free of the bandwidth. Then a global R-squared for VCM is defined as

(3.6) $$R_p^2(h) = 1 - \frac{SSE_p(h)}{SST} = \frac{SSR_p(h)}{SST}.$$

To investigate the asymptotic properties of the global ANOVA R-squared (3.6), we impose Conditions (A2), (A3), (A5), and the following technical conditions:

CONDITIONS (B).

- (B1) The second derivative of the density $g(u)$ is bounded, continuous, and square integrable on a compact support.
- (B2) The $(p+1)$st derivative of $a_j(\cdot)$, $j = 1, \ldots, d$, exists.
- (B3) $EX_j^{2s} < \infty$, for some $s > 2$, $j = 1, \ldots, p$.
- (B4) Let $\gamma_{ij}(u) = E(X_i X_j | U = u)$, $i,j = 1, \ldots, d$, $\gamma_{ij}(\cdot)$ is continuous in a neighborhood of $u$.

Now, we state the asymptotic normality for the global ANOVA R-squared (3.6) in the following theorem and its proof is given in the Appendix.

THEOREM 6. *Assume that as $n \to \infty$, $h = h(n) \to 0$. When fitting LPR with an odd $p$, under Conditions* (A2), (A3), (A5) *and* (B1)–(B4), *with $nh^{2p+2} \to 0$ and $nh^2 \to \infty$:*

(a) *The asymptotic conditional bias of $R_p^2(h)$ is*

$$-h^2 \frac{\mu_2}{2\sigma_y^2} \int \sigma^2(u) g''(u)\, du (1 + o_P(1)).$$

(b) *The asymptotic conditional variance of $R_p^2(h)$ is*

$$n^{-1}\left(\frac{\mathrm{Var}(\varepsilon^2)}{\sigma_y^4} E(\sigma^4(U))\left(\int K_0^*(v)\, dv\right) + \frac{(m_4 - \sigma_y^4)(E(\sigma^2(U)))^2}{\sigma_y^8}\right)$$
$$\times (1 + o_P(1)).$$

(c) *Under the homoscedastic assumption and conditioned on $\{X_1, \ldots, X_n\}$, $R_p^2(h)$ converges in distribution to a normal distribution with the above asymptotic conditional bias and variance.*

(d) *Under the assumptions in* (c), *$SSE_p(h)$ for VCM is a $\sqrt{n}$-consistent estimate for $\sigma^2$. Its asymptotic conditional bias is $o_P(n^{-1/2})$ if $\int g''(u)\, du = 0$ and the asymptotic conditional variance $n^{-1}\sigma^4(\int K_0^*(v)\, dv)(1 + o_P(1))$.*

Theorem 6 extends Theorem 2 to VCM. Other ANOVA results for VCM, such as degree of freedom, testing against $H_0 : a_k(U) = c$ for some $k$ with $c$ a constant, and testing for overall model significance, will be derived in a separate paper.

**4. Numerical results.** In this section, we use computer simulations to investigate the performance of the ANOVA R-squared and the proposed $F$-test.

4.1. *Simulation results for the ANOVA R-squared.* Two examples from Doksum and Froda [6] are used to compare the performance between the ANOVA R-squared (2.10), the adjusted ANOVA R-squared (2.18), the correlation R-squared (1.3), and an empirical RSS-related R-squared $R_s^2 = RSS/\sum_i (Y_i - \bar{Y})^2$. For comparison only, we also include the R-squared from fitting a simple linear model. Sample sizes of $n = 50$ and 200 are used with 400 simulations. Following Doksum and Froda [6], we use a fixed bandwidth $h = 0.22$ (approximately 0.7 times the standard deviation of $X$ in the examples). The purpose is to see how the four coefficients of determination differ from one another when the same amount of smoothing is applied. Local linear regression with the Epanechnikov kernel $K(u) = 0.75(1-u^2)I_{|u|\leq 1}$ is applied and 200 equally spaced grid points on $(\min_i X_i, \max_i X_i)$ are used to approximate the integration for $R_1^2(h)$ and $R_{1,adj}^2(h)$. No special treatment for boundary points is implemented for any of the four nonparametric R-squared's.

EXAMPLE 1. Bump model: $Y = 2 - 5(X - e^{-100(X-0.5)^2}) + \sigma\varepsilon$, where $X$ follows a Uniform$(0,1)$ and the distribution of $\varepsilon$ is $N(0,1)$. $X$ and $\varepsilon$ are independent, and $\sigma = 0.5, 1, 2, 4$ results in high to low values for the true value of the coefficient of determination. The results show that the four nonparametric R-squared's have similar performance for both $n=50$ and $n=200$, and hence the plots are omitted for brevity. The values for the ANOVA R-squared is slightly smaller than $R_\rho^2$ and $R_s^2$; for example, when $\sigma = 0.5$, the average $R_1^2$ is 0.8155 (sd 0.0325), 0.8273 (sd 0.0323) for $R_\rho^2$, and 0.8337 (sd 0.0337) for $R_s^2$.

EXAMPLE 2. Twisted pear model: $Y = 5 + 0.1Xe^{(5-0.5X)} + \frac{(1+0.5X)}{3}\sigma\varepsilon$, where $X \sim N(1.2, (1/3)^2)$ and $\varepsilon \sim N(0,1)$. $X$ and $\varepsilon$ are independent, and the values of $\sigma$ are the same as in Example 1. The original model from Doksum and Froda [6] did not include the constant 5. We add a nonzero constant in the model for convenience of performing $F$-tests in Section 4.2. This model represents a situation where the relationship between $X$ and $Y$ is strong for small $x$, but then tapers off as the noise variance increases. Figure 1 gives the boxplots for $n = 50$. Clearly both the unadjusted and adjusted ANOVA R-squared's behave much more stably than $R_\rho^2$ and $R_s^2$. When $\sigma = 0.5$, the values of mean (sd) are 0.9512 (0.0195), 0.9444 (0.0216), 0.8587 (0.1662) and 0.8730 (0.1752) for $R_1^2$, adjusted $R_{1,\text{adj}}^2$, $R_\rho^2$ and $R_s^2$, respectively. Both $R_\rho^2$ and $R_s^2$ have a skewed distribution for this heteroscedastic model. Similar results can be observed for the case of $\sigma = 1$. When $\sigma = 4$, we note that there is one negative $R_s^2$ and four negative $R_{1,\text{adj}}^2$, which are not guaranteed to lie between 0 and 1. The results for $n = 200$ are similar to those of $n = 50$ and hence are omitted. This example demonstrates some advantages of the ANOVA R-squared in a heteroscedastic model as compared to other nonparametric coefficients of determination.

4.2. *Simulation results for the F-test of no effect.* Due to boundary effects in practice, we adopt a more conservative version of the $F$-statistic, defined as

$$(4.1) \qquad F(h) = \frac{(SSR_p(h)/(\text{tr}(H^*)-1)}{(\sum_i (Y_i - \bar{Y})^2 - SSR_p(h))/(n - \text{tr}(H^*))},$$

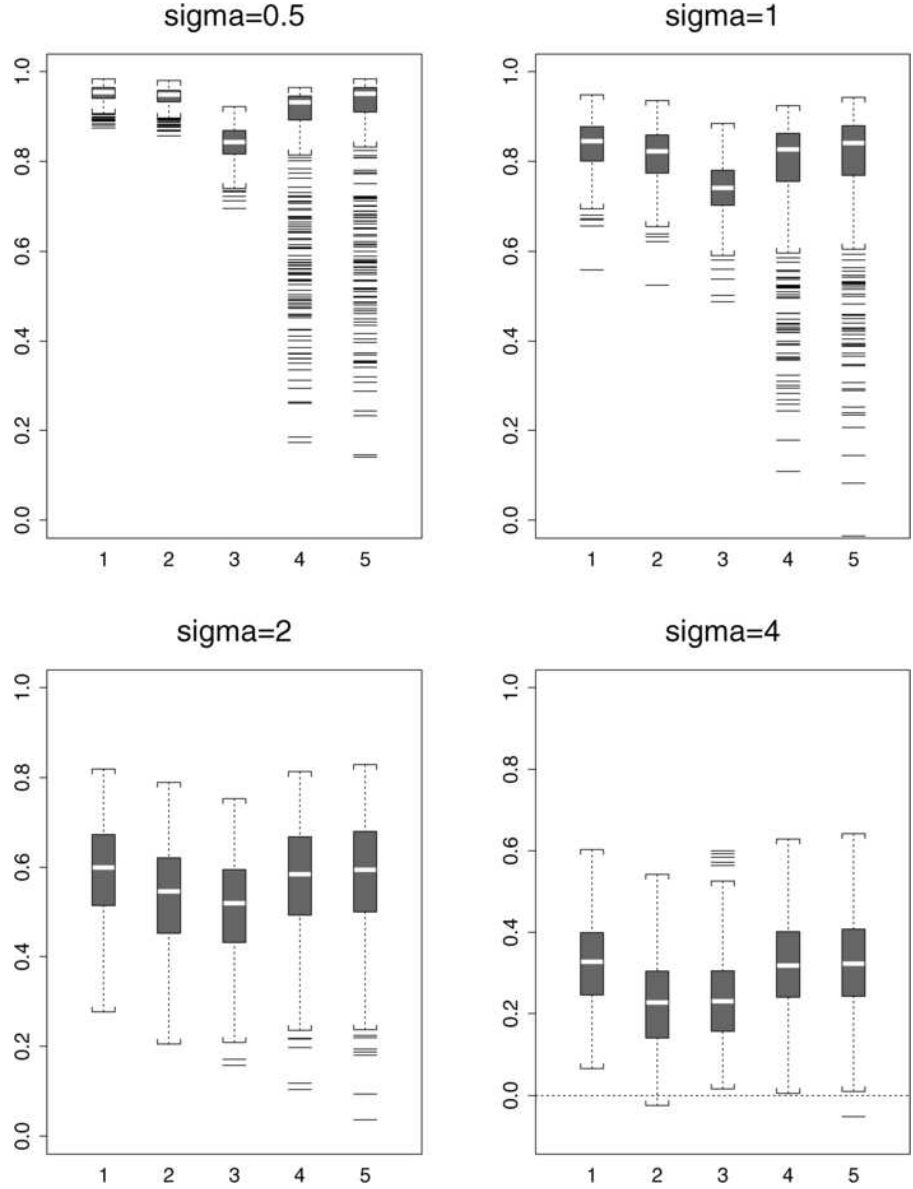

FIG. 1. *Example 2. Boxplots for 400 trials of five different R-squared's with $n = 50$*: (1) *ANOVA* $R_1^2$, (2) *adjusted ANOVA* $R_1^2$, (3) *R-squared from fitting a simple linear model*, (4) $R_\rho^2$ *(Doksum and Samarov [5]) and* (5) *empirical* $R_s^2$.

where $SSR_p(h)$ is estimated based on (2.8) without any boundary adjustment. Note that in the denominator of (4.1), $(\sum_i (Y_i - \bar{Y})^2 - SSR_p(h))$ is used instead of $SSE_p(h)$. Examples 1 and 2 with $\sigma = 1$ are modified as Examples 3 and 4 to illustrate the proposed $F$-test. For each example, three fixed values of the bandwidth are used: Example 3, $h = 0.15$, 0.22 and 0.34, and Example 4, $h = 0.22$, 0.34 and 0.51, with a ratio of roughly 1.5. The $F$-test statistic in (4.1) is calculated and its $p$-value is obtained using the $F$-distribution with degrees of freedom $(\text{tr}(H^*) - 1, n - \text{tr}(H^*))$. A significance level 0.05 is used to determine whether to reject the null hypothesis or not. Again sample sizes $n = 50$ and $n = 200$ are used with 400 simulations. For comparison only, we also include another $F$-flavor test, the pseudo-likelihood ratio test (PLRT) for no effect by Azzalini and Bowman [2], in which a chi-squared distribution was calibrated to obtain the $p$-value.

EXAMPLE 3. Consider the model: $Y = 2 - a \times (X - e^{-100(X-0.5)^2}) + \varepsilon$, where $a = 0, 0.5, \ldots, 3$, $X \sim \text{Uniform}(0,1)$ and $\varepsilon \sim N(0,1)$. The case of $a = 0$ gives the intercept only model, while Example 1 corresponds to $a = 5$. Figure 2(a) illustrates the shapes of the true regression functions. The proportions of rejection by the $F$-statistic (4.1) and PLRT are plotted in Figure 2(b)–(d) as a function of $a$. With a conservative (4.1), all type-I errors of the $F$-test are below 5% level. The PLRT has slightly better power than the $F$-test when $n = 50$, and the two tests behave similarly for $n = 200$. Both tests have better power with bandwidth increasing, while the type-I error of the PLRT exceeds 0.05 level when $h = 0.34$.

EXAMPLE 4. Consider the following model: $Y = 5 + aXe^{(5-0.5X)} + \frac{(1+0.5X)}{3}\varepsilon$, where $a = 0, 0.01, \ldots, 0.06$, $X \sim N(1.2(1/3)^2)$, and $\varepsilon \sim N(0,1^2)$. For this heteroscedastic model, $a = 0$ corresponds to the null hypothesis, and Example 2 corresponds to $a = 0.1$. We note that neither of the two tests is formally applicable, but we want to examine their robustness against deviations from homoscedasticity. A plot of the true regression functions is given in Figure 2(e), and the percentages of rejection over 400 simulations are given in Figure 2(f)–(h). As in Example 3, the PLRT has slightly better power than the $F$-test when $n = 50$. We observe a less accurate approximation of the type-I error by the PLRT when $n = 200$: 7.75%, 6.5% and 6.25% for $h = 0.22$, 0.34 and 0.51, respectively (the corresponding numbers are 4.5%, 4% and 4% for the $F$-test). This may justify PLRT's better performance when $a = 0.01$ and 0.02. This example shows that even under a heteroscedastic error structure, both tests perform reasonably well.

**5. Real data.** The data from Simonoff [25] were obtained in Lake Erie, containing 52 rows numbered consecutively from the northwest (row 1) to the southeast (row 52) and the sum of yields of the harvest in 1989, 1990 and 1991, as measured by the total number of lugs (a lug is roughly 30 pounds of grapes). Figure 3(a) shows the data and the local linear estimates at grid points $1, 1.5, \ldots, 52$, with the Gaussian kernel and bandwidth $h = 3$ (solid line) and $h = 1.5$ (dashed line). The choice of bandwidth follows Simonoff [25]. The dip in yield around rows 30–40 is possibly due to a farmhouse directly opposite those rows (Simonoff [25]). The coefficients of determination indicate good explanatory power: when $h = 3$, $R_1^2$, $R_\rho^2$ and $R_s^2$ are 0.8493, 0.9414, 0.8854, respectively; when $h = 1.5$, 0.9046, 0.9638 and 0.9297. The corresponding pointwise R-squared is shown in Figure 3(b). The curve with $h = 1.5$ has a larger pointwise $R_1^2(x)$ in most locations than that of $h = 3$. The local R-squared with $h = 3$ is only 40–50% for rows 31–34, and above 90% for rows 12–23 and 46–52, reflecting some difference across data in the proportion of variation explained by local linear regression. The difference leads to the idea of using the local R-squared for variable bandwidth selection in a future paper. The ANOVA tables for $h = 3$ and 1.5 are given in Tables 2 and 3. As expected, the $SSR_1$ of $h = 1.5$ is greater than that of $h = 3$. Both $p$-values of the ANOVA $F$-statistic (4.1) are $<10^{-7}$, indicating rejection of the null hypothesis. The PLRT also gives very small $p$-values, $4.3 \times 10^{-4}$ and $1.33 \times 10^{-4}$ for $h = 1.5$ and 3, respectively. Note that due to boundary effects, $SSR_p(h) + SSE_p(h)$ does not equal the sample variance of $Y$. We give both quantities in the ANOVA tables to illustrate this effect in practice.

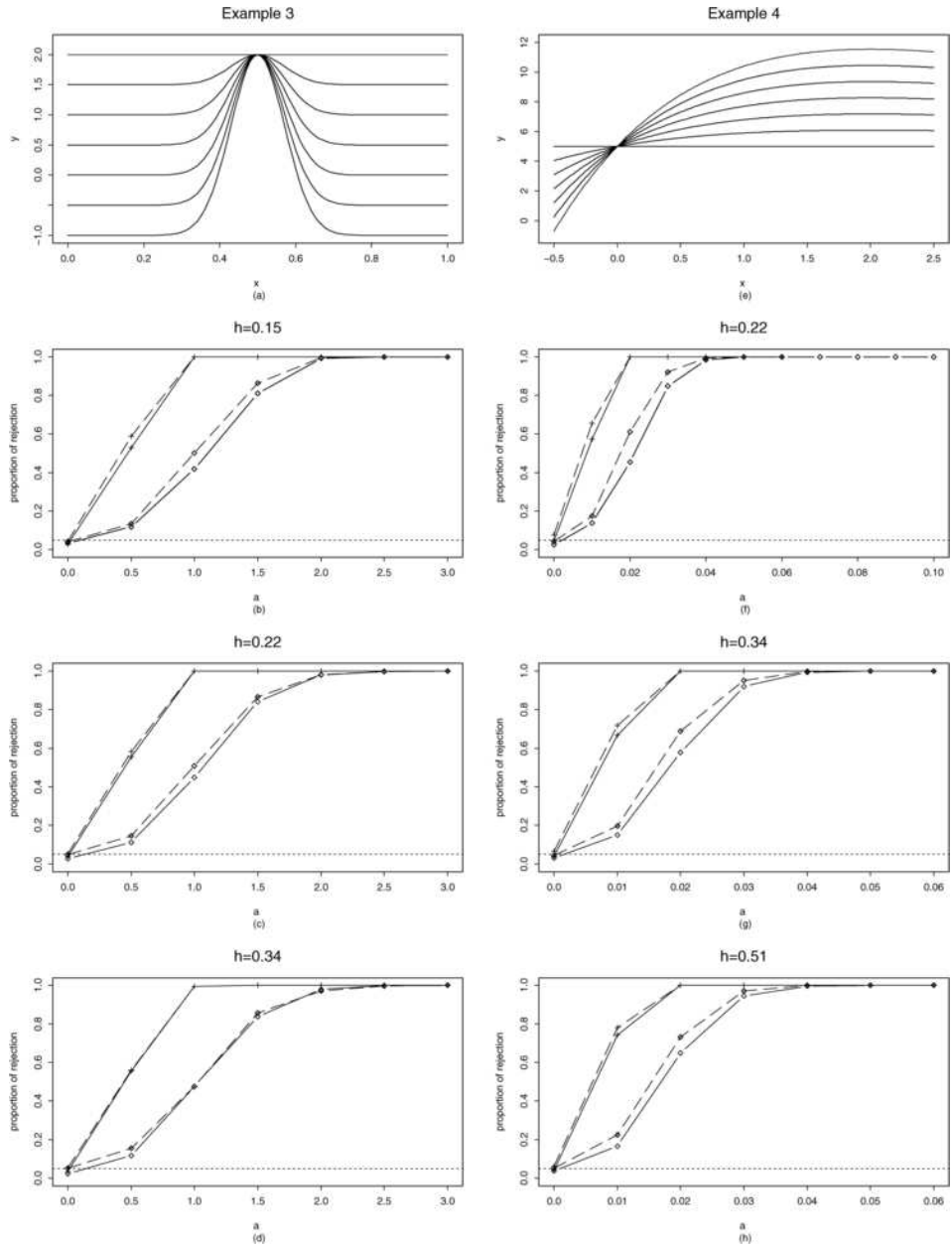

Fig. 2. *Examples 3 and 4.* (a) *Plot of the true regression function curves for Example 3, $a = 0, 0.5, \ldots, 3$;* (b) *the percentages of rejection for simulated data in Example 3 based on 400 simulations with $h = 0.15$; ANOVA F-test (solid line), PLRT (long dash line), and the short dash line indicates the 5% significance level;* $+$ *($n = 200$);* $\circ$ *($n = 50$);* (c) *same as in* (b) *except $h = 0.22$;* (d) *same as in* (b) *except $h = 0.34$;* (e) *plot of the true regression function curves for Example 4, $a = 0, 0.01, \ldots, 0.06$;* (f)–(h): *same as in* (b)–(d) *for Example 4 with $h = 0.22$, 0.34 and 0.51, respectively.*

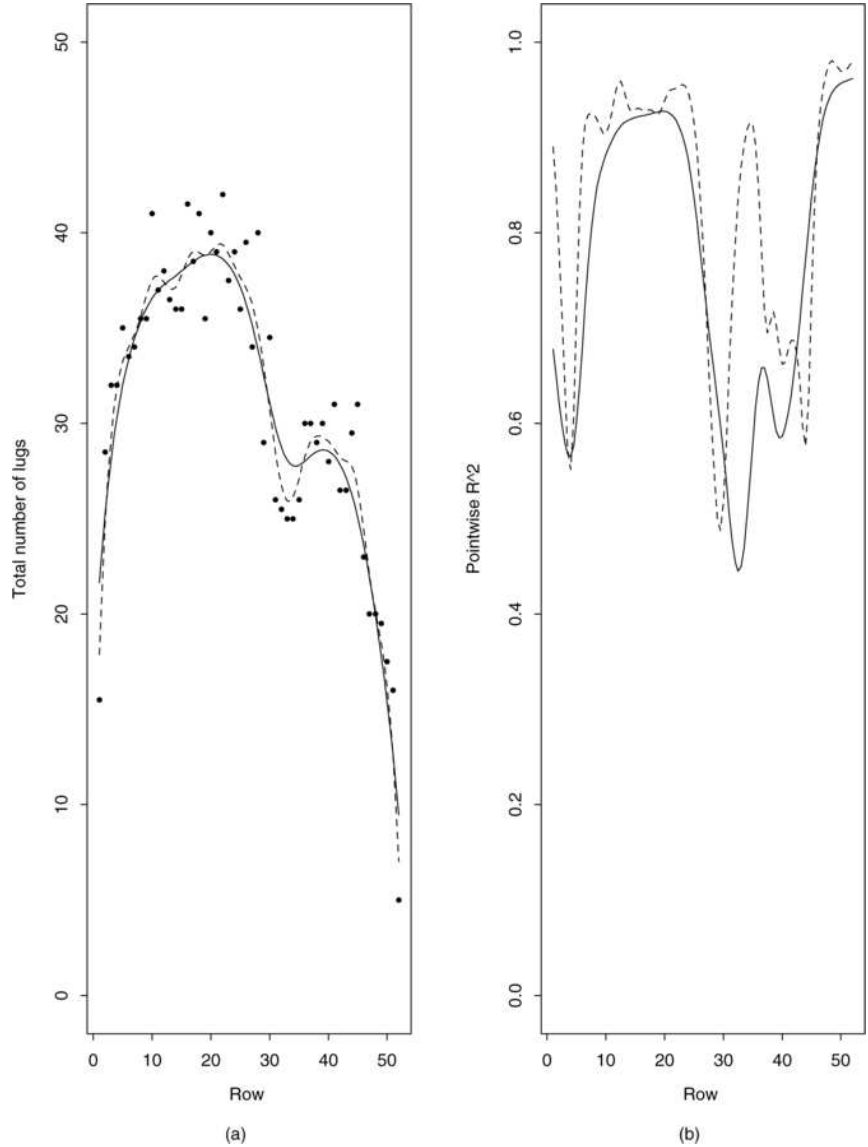

FIG. 3. (a) *Scatterplot of total lug counts versus row for the vineyard data with local linear estimates $h = 3$ (solid line) and $h = 1.5$ (dashed line).* (b) *Plot of the corresponding pointwise $R_1^2(x)$.*

TABLE 2
*ANOVA table for vineyard data with bandwidth 3*

| Source | Degree of freedom | Sum of squares | Mean squares | $F$ |
|---|---|---|---|---|
| Regression | $(7.5509 - 1)$ | $SSR_p = \frac{2204.6682}{52}$ | $\frac{2204.6682}{6.5509}$ | $F = 15.3299$ |
| Residual | $(52 - 7.5509)$ | $SSE_p = \frac{391.2489}{52}$ | $\frac{391.2489}{44.4491}$ | |
| Total | 51 | $SST = \frac{2595.9171}{52}$ | $\frac{3180.4808}{52}$ | |

**6. Discussion.** Though the idea of nonparametric ANOVA inference is not new, we believe that the work in this paper provides a unified framework with an asymptotic geometric configuration for the first time. The proposed ANOVA tools for LPR are easy to carry out in practice and we hope that the methodology will be useful for data analysis. It will be interesting to explore a similar ANOVA framework for other nonparametric regression methods such as penalized splines in future studies. The ground-breaking points are the elegant local ANOVA decomposition (2.5) and construction of global ANOVA quantities through integrating local counterparts. Thus LPR, fitting local polynomials across data, may be viewed as a "calculus" extension of classical polynomial models. A surprising by-product is that the projected response $H^*\mathbf{y}$ has a bias of order $h^4$, which is smaller than the usual order $h^2$. The proposed projection matrix $H^*$ overcomes the problem of a nonsymmetric smoother matrix of local linear regression, and we show that it has nice geometric properties that lead to a natural $F$-test for no-effect. $H^*$ also provides a new geometric view of LPR: for example, in the case of local linear regression, the local fitting at $x$ is to project $\mathbf{y}$ into local column space of $\mathbf{X}$ and the locally projected values are $\hat{\beta}_0(x) + \hat{\beta}_1(x)(X_i - x)$, $i = 1, \ldots, n$; these locally projected values at different grid points around $X_i$ are then combined through weighted integration to form the projected value $Y_i^*$ [see (2.15)]. The projection view and the geometric representation of the ANOVA quantities offer new insights for LPR. The proposed $F$-test shares the property of the "Wilks phenomenon" with the generalized likelihood ratio test (Fan, Zhang and Zhang [10]), in that it does not depend on nuisance parameters. The numerical results presented in the paper show that the test statistic under the null hypothesis follows well the asymptotic $F$-distribution without further calibration; one does not have to simulate the null distributions to obtain the critical value. The paper also presents a brief multivariate extension of nonparametric ANOVA inference to VCM; more details will be developed in a separate paper. Based on findings in this paper, several follow-up problems are being investigated, including extension of the $F$-test to test for a polynomial relationship (Huang and Su [16]), and ANOVA inference for partial linear models and generalized additive models. We are also interested in applying the ANOVA approach to study the bandwidth selection problem, for example, using the local R-squared for variable bandwidth selection, and using the classical model selection criteria of AIC and BIC with the proposed $SSE_p(h)$ and degree of freedom $\text{tr}(H^*)$ for global bandwidth selection.

TABLE 3
*ANOVA table for vineyard data with bandwidth 1.5*

| Source | Degree of freedom | Sum of squares | Mean squares | $F$ |
|---|---|---|---|---|
| Regression | $(14.8095 - 1)$ | $SSR_p = \frac{2461.9695}{52}$ | $\frac{2461.9695}{13.8095}$ | $F = 8.9798$ |
| Residual | $(52 - 14.8095)$ | $SSE_p = \frac{259.5325}{52}$ | $\frac{259.5325}{36.1905}$ | |
| Total | 51 | $SST = \frac{2721.5020}{52}$ | $\frac{3180.4808}{52}$ | |

## APPENDIX

Proofs of Theorems 3, 4 and 6 are included in this section. The following lemma by Mack and Silverman [18] will be needed.

LEMMA A.1. *Assume that $E|Y^3| < \infty$ and $\sup_x \int |y|^s f(x,y)\,dy < \infty$, where $f(x,y)$ denotes the joint density of $(X,Y)$. Let $K$ be a bounded positive function with a bounded support, satisfying a Lipschitz condition, and $D$ the support for the marginal density of $X$. Then*

$$\sup_{x \in D}\left|n^{-1}\sum_{i=1}^n \{K_h(X_i - x)Y_i - E[K_h(X_i - x)Y_i]\}\right|$$
$$= O_P[\{nh/\log(1/h)\}^{-1/2}],$$

*provided that $n^{2a-1}h \to \infty$ for some $a < 1 - s^{-1}$.*

PROOF OF THEOREM 3. For the $i$th element $Y_i^*$, under Conditions (A1) and (A2),

$$
\begin{aligned}
Y_i^* - m(X_i) &= \int \left(\sum_{j=0}^p (\hat{\beta}_j(x)(X_i - x)^j)\right) K_h(X_i - x)\,dx \\
&\quad - \int m(X_i) K_h(X_i - x)\,dx \\
&= \int ((\hat{\beta}_0(x) - \beta_0(x)) + \cdots + (\hat{\beta}_p(x) - \beta_p(x))(X_i - x)^p) \\
&\quad \times K_h(X_i - x)\,dx \\
&\quad - \int (\beta_{p+1}(x)(X_i - x)^{p+1} + r(x, X_i)) K_h(X_i - x)\,dx,
\end{aligned}
\tag{A.1}
$$

where $r(x, X_i)$ denotes the remainder terms after a $(p+1)$st-order Taylor expansion. By using the bias expression from Wand and Jones [26], for example, when $p$ is odd,

$$E\left\{\int (\hat{\beta}_0(x) - \beta_0(x)) K_h(X_i - x)\,dx \Big| X_1, \ldots, X_n\right\} = h^{p+1} b_{0,p}(x)(1 + o_P(1)),$$

and similarly for $\int (\hat{\beta}_j(x) - \beta_j(x))(X_i - x)^j K_h(X_i - x)\,dx$, $j \geq 1$. With

$$
\begin{aligned}
&\int \beta_{p+1}(x)(X_i - x)^{p+1} K_h(X_i - x)\,dx \\
&\quad = \frac{1}{(p+1)!} h^{p+1} \mu_{p+1} m^{(p+1)}(X_i)(1 + o_P(1)),
\end{aligned}
$$

the asymptotic conditional bias of $Y_i^*$ in (2.12) is obtained when $p$ is odd. The case for an even $p$ follows analogously. For local linear regression in part (b), the $h^2$-order terms are canceled, and the asymptotic conditional bias follows from further expansion of (A.1).

For part (c), denote the conditional bias vector of $\mathbf{y}^*$ by $\mathbf{b} = H^*\mathbf{m} - \mathbf{m}$. Then

$$E\{H^{*^2}\mathbf{y}|X_1, \ldots, X_n\} - \mathbf{m} = H^*(\mathbf{m} + \mathbf{b}) - \mathbf{m} = \mathbf{b} + H^*\mathbf{b}. \tag{A.2}$$

The rate of $\mathbf{b} = H^*\mathbf{m} - \mathbf{m}$ is given in (2.12) and (2.13). It remains to investigate the rate of elements in $H^* = (h_{i,j}^*)$. The $(j,k)$th element of $(\mathbf{X}^T W \mathbf{X})$ matrix is $s_{j,k}(x) = \sum_i (X_i - x)^{j+k-2} K_h(X_i - x)/\hat{f}(x) = nh^{j+k-2}(\mu_{j+k-2} + \mu_{j+k-1} f'(x)/f(x))(1 + o_P(1))$ by Lemma A.1, and

$$\mathbf{X}^T W \mathbf{X} = nD(S_p + hS_p' f'(x)/f(x) + o_P(h^2))D, \tag{A.3}$$

where $D = \text{diag}(1, h, \ldots, h^p)$, $S_p = (\mu_{i+j-2})_{1 \leq i,j \leq (p+1)}$ and $S'_p = (\mu_{i+j-1})_{1 \leq i,j \leq (p+1)}$. Then $(\mathbf{X}^T W \mathbf{X})^{-1} = D^{-1} S_p^{-1} D^{-1}(1 + o_P(1))$. Denote $S_p^{-1} = (s_{i,j})$ and $s_{i,j}$ is of order $O(1)$. Then before integrating over $x$, the $i$th diagonal element of $WH\hat{f}(x) = W\mathbf{X}(\mathbf{X}^T W \mathbf{X})^{-1}\mathbf{X}^T W \hat{f}(x)$ has a form:

$$n^{-1} \sum_{l=0}^{p} \sum_{k=0}^{p} (K_h(X_i - x))^2 (X_i - x)^{l+k} s_{(k+1),l} h^{-(l+k)}(1 + o_P(1)),$$

which is of order $O(n^{-1}h^{-1})(1 + o_P(1))$. After integration, the rate for $h^*_{i,i}$ remains $O(n^{-1}h^{-1})(1 + o_P(1))$. We next show that the rate of nondiagonal elements of $H^*$ is of order $O(n^{-1})(1 + o_P(1))$. For $i \neq j$, the integrand for $h^*_{i,j}$ is

$$\sum_{l=0}^{p} \sum_{k=0}^{p} K_h(X_i - x) K_h(X_j - x)(X_i - x)^l (X_j - x)^k s_{(k+1),l} h^{-(l+k)},$$

which is of order $O(n^{-1})(1 + o_P(1))$. Then results stated in (c) follow from (A.2).

For part (d), under the homoscedastic model, $\text{Var}(\mathbf{y}^*) = H^{*2}\sigma^2$, and the conditional variance of $Y_i^*$ is $\sigma^2 \sum_j h_{i,j}^{*2}$. When $i = j$, $h_{i,i}^{*2}$ is of order $O(n^{-2}h^{-2})(1 + o_P(1))$. For $i \neq j$,

$$\sum_{j \neq i} h_{i,j}^{*2} = n^{-1}h^{-2} \int \left\{ \int K(u) K\left(\frac{X_j - X_i}{h} - u\right)\left(1 - \frac{u}{\mu_2}\left(\frac{X_j - X_i}{h} - u\right)\right) du \right\}^2$$

$$\times f(X_j) \, dX_j (1 + o_P(1))$$

$$= n^{-1}h^{-1} \frac{1}{f(X_i)} \int \left\{ K_0^*(u) - \frac{K_1^*(u)}{\mu_2} \right\}^2 du (1 + o_P(1)).$$

Hence the asymptotic conditional variance of $Y_i^*$ is as given in (2.14). This completes the proof of Theorem 3. $\square$

PROOF OF THEOREM 4. Theorem 4(a) follows directly from Theorem 3. For part (b), since $H^*$ is asymptotically idempotent, $(H^* - L)$ is asymptotically idempotent given that $(H^* - L)^2 = H^{*2} - L$. Therefore with the homoscedastic normality assumption under the no-effect null hypothesis, $SSR_p(h)$ has an asymptotic $\chi^2$-distribution with degree of freedom $(\text{tr}(H^*) - 1)$. Similarly for $SSE_p(h)$, it has an asymptotic $\chi^2$-distribution with a degree of freedom $(n - \text{tr}(H^*))$. With $(H^* - L)$ and $(I - H^*)$ being asymptotic orthogonal in part (a), the test statistic $F$ in (2.16) has an asymptotic $F$-distribution with a degree of freedom $(\text{tr}(H^*) - 1, n - \text{tr}(H^*))$.

For part (c), note that $\text{tr}(H^T W \hat{f}(x)) = \text{tr}(\hat{f}(x) W \mathbf{X}(\mathbf{X}^T W \mathbf{X})^{-1} \mathbf{X}^T W) = \text{tr}((X^T W X)^{-1} X^T W^2 X \hat{f}(x))$. Using $S_p$ in (A.3) with $p = 1$, and

$$(\hat{f}(x) \mathbf{X}^T W^2 \mathbf{X}) = \begin{pmatrix} \nu_0/h & \nu_2 f'(x)/f(x) \\ \nu_2 f'(x)/f(x) & h\nu_2 \end{pmatrix}(1 + o_P(1)),$$

(2.17) is obtained. Therefore the proof of Theorem 4 is complete. □

PROOF OF THEOREM 6. We need the following notation: the error vector as $\mathbf{e} = (\sigma(U_1)\varepsilon_1, \ldots, \sigma(U_n)\varepsilon_n)^T$, the mean vector as $\mathbf{m} = (\mathbf{a}(U_1)^T\mathbf{X}_1, \ldots, \mathbf{a}(U_n)^T\mathbf{X}_n)^T$ with $\mathbf{X}_i = (1, X_{i2}, \ldots, X_{i,d})^T$, and $D = I_{p+1} \otimes \text{diag}(1, \ldots, h^p)$. Let $\mu = (\mu_{p+1}, \ldots, \mu_{2p+1})^T$ [recall $\mu_j = \int u^j K(u)\,du$]. It follows from (3.6) that

$$(A.4) \quad (R_p^2(h) - \eta) = \frac{1}{SST}\left\{-(SSE_p(h) - E(\sigma^2(U))) + (SST - \sigma_y^2)\frac{E(\sigma^2(U))}{\sigma_y^2}\right\}.$$

The first term in (A.4) can be expressed as

$$SSE_p(h) - E(\sigma^2(U)) = \frac{1}{n}\int \mathbf{y}^T(W_u - W_u H_u)\mathbf{y}\hat{g}(u)\,du - E(\sigma^2(U))$$

$$\equiv I_1 + I_2 + I_3,$$

where $I_1 = \frac{1}{n}\int \mathbf{e}^T(W_u - W_u H_u)\mathbf{e}\hat{g}(u)\,du - E(\sigma^2(U))$, $I_2 = \frac{1}{n}\int \mathbf{m}^T(W_u - W_u H_u)\mathbf{m}\hat{g}(u)\,du$ and $I_3 = \frac{2}{n}\int \mathbf{e}^T(W_u - W_u H_u)\mathbf{m}\hat{g}(u)\,du$.

For matrix $H_u$, using Lemma A.1 we find that $\mathbf{X}_u^T W_u \mathbf{X}_u = D(\Gamma \otimes S_p)D(1 + o_P(1))$ and $\sum_{k=1}^d a_k^{(p+1)}(u)\mathbf{X}_u^T W_u(X_{1k}(U_1-u)^{p+1}, \ldots, X_{nk}(U_n-u)^{p+1})^T = D(\Gamma \otimes \mu)\mathbf{a}^{(p+1)T}h^{p+1}(1 + o_P(1))$, where $\Gamma = E\{(X_1, \ldots, X_d)^T(X_1, \ldots, X_d) | U = u\}$ and $\mathbf{a}^{(p+1)} = (a_1^{(p+1)}, \ldots, a_d^{(p+1)})^T$. The term $I_2$ conditioned on $\{X_1, \ldots, X_n\}$ is nonrandom and asymptotically $I_2 = \frac{1}{\{(p+1)!\}^2}h^{2p+2}g(u)(\mu_{2p+2} - \mu^T S_p^T \times \mu_{2p+2})(1 + o_P(1))$. Conditioned on $\{X_1, \ldots, X_n\}$, $I_3$ has a mean 0 and its variance is of order $h^{2(p+1)}$. For $I_1$, assuming local homoscedasticity,

$$(A.5) \quad I_1 = \frac{1}{n}\sum_{i=1}^n\left[\int K(U_i - u)\sigma^2(u)\varepsilon_i^2\,du\right] - E(\sigma^2(U)) + O_P(n^{-1}h^{-1}).$$

The asymptotic conditional mean and variance for $I_1$ are

$$E(I_1|X_1, \ldots, X_n) = h^2\frac{\mu_2}{2}\int \sigma^2(u)g''(u)\,du(1 + o_P(1)),$$

$$(A.6) \quad \text{Var}(I_1|X_1, \ldots, X_n) = n^{-2}\sum_i\left(\int \sigma^2(u)K_h(U_i - u)\,du\right)^2$$

$$= n^{-1}E(\sigma^4(U))\left(\int K_0^*(v)\,dv\right)(1 + o_P(1)).$$

It is clear that under the condition $nh^{2p+2} \to 0$, $I_1$ is the dominating term for $(SSE_p(h) - E(\sigma^2(U)))$. Further, the asymptotic conditional variance of $I_1$ is dominated by (A.6) since $O_P(n^{-1}h^{-1})$ in (A.5) is smaller than (A.6) under the condition that $nh^2 \to \infty$. Using Theorem 8.16 in Lehmann and Casella [17], $\sqrt{n}(SST - \sigma_y^2)$ has the asymptotic normality $N(0, \text{Var}[(Y - \mu_y)^2])$. Then from (A.4), the asymptotic conditional variance of $(R_p^2(h) - \eta^2)$ is obtained.

Last we establish the asymptotic normality of $R_p^2(h)$ in the homoscedastic case. Since $SST \to \sigma_y^2$ with probability 1, (A.4) becomes

$$(R_p^2(h) - \eta) = \left\{-\frac{1}{\sigma_y^2} I_1(1 + o_P(1)) + (SST(h) - \sigma_y^2)\frac{1}{\sigma_y^4}\right\}.$$

$I_1$ is a summation of i.i.d. random variables $\varepsilon_i^2$, $i = 1, \ldots, n$, and hence by the central limit theorem, $I_1$ has an asymptotic normal distribution. It is easy to show that the covariance of $I_1$ and $SST$ conditioned on $X_1, \ldots, X_n$ is of smaller order than the sum of variances of $I_1$ and $(SST - \sigma_y^2)$. Thus, the asymptotic normality for $R_p(h)$ is obtained. The results in part (d) are easily seen from asymptotic normality of $I_1$. □

**Acknowledgments.** The authors are grateful to two referees, an Associate Editor and Editors for helpful comments. Thanks are due also to Dr. Hua He for comments on the first version.

DEPARTMENT OF BIOSTATISTICS
AND COMPUTATIONAL BIOLOGY
UNIVERSITY OF ROCHESTER
601 ELMWOOD AVENUE
BOX 630
ROCHESTER, NEW YORK 14642
USA
E-MAIL: lhuang@bst.rochester.edu

DEPARTMENT OF MATHEMATICS AND STATISTICS
SAN DIEGO STATE UNIVERSITY
5500 CAMPANILE DRIVE
SAN DIEGO, CALIFORNIA 92182
USA
E-MAIL: jchen@sciences.sdsu.edu